\documentclass[letterpaper, 10 pt, conference]{ieeeconf}  

\IEEEoverridecommandlockouts                              
\title{\LARGE \bf
A Consensus-ADMM Approach for Strategic Generation Investment in Electricity Markets
}

\author{Vladimir Dvorkin, Jalal Kazempour, Luis Baringo and Pierre Pinson
\thanks{The work of L. Baringo has been partially funded
by the Ministry of Science of Spain under CICYT Project
ENE2015-63879-R (MINECO/FEDER,UE).}
\thanks{V. Dvorkin, J. Kazempour. P. Pinson are with the Department of Electrical Engineering,
        Technical University of Denmark, Lyngby, Denmark
        {\tt\small \{vladvo,seykaz,ppin\}@elektro.dtu.dk}}%
\thanks{L. Baringo is with the School of Industrial Engineering, Universidad de Castilla - La Mancha, Ciudad Real, Spain
        {\tt\small luis.baringo@uclm.es}}%
}

\usepackage{graphicx}
\usepackage{mathtools}
\usepackage{amssymb}
\usepackage{amsfonts}
\usepackage{amsmath}
\usepackage{nccmath}
\usepackage{bm}
\usepackage{soul}
\usepackage{multirow}
\usepackage[ruled]{algorithm2e}
\usepackage{ragged2e}
\usepackage{makecell}
\usepackage{tikz}
\usetikzlibrary{shapes.geometric}
\usetikzlibrary{decorations.pathreplacing}
\usepackage{pgfplots}
\usepackage{url}
\usepackage[flushleft]{threeparttable}
\usepackage{verbatim}

\newenvironment{ldescription}[1]
{\begin{list}{}%
{\renewcommand\makelabel[1]{##1\hfill}%
\settowidth\labelwidth{\makelabel{#1}}%
\setlength\leftmargin{\labelwidth}
\addtolength\leftmargin{\labelsep}}}
{\end{list}}

\newcommand\norm[1]{\left\lVert#1\right\rVert}

\newcolumntype{M}[1]{>{\centering\arraybackslash}m{#1}}
\newcolumntype{N}{@{}m{0pt}@{}}

\definecolor{color_ConvD}{rgb}{0.8588,0.3333,0.3333}
\definecolor{color_bilevel}{rgb}{0.3216,0.6863,0.8980}
\definecolor{color_StochD}{rgb}{0.333,0.8588,0.6039}

\definecolor{blue}{rgb}{0,0,0}

\begin{document}

\maketitle
\thispagestyle{empty}
\pagestyle{empty}

\begin{abstract}
This paper addresses a multi-stage generation investment problem for a strategic (price-maker) power producer in electricity markets. This problem is exposed to different sources of uncertainty, including short-term operational (e.g., rivals' offering strategies) and long-term macro (e.g., demand growth) uncertainties. This problem is formulated as a stochastic bilevel optimization problem, which eventually recasts as a large-scale stochastic mixed-integer linear programming (MILP) problem with limited computational tractability. To cope with computational issues, we propose a consensus version of alternating direction method of multipliers (ADMM), which decomposes the original problem by both short- and long-term scenarios. Although the convergence of ADMM to the global solution cannot be generally guaranteed for MILP problems, we introduce two bounds on the optimal solution, \textcolor{blue}{allowing for the evaluation of the solution quality over iterations}. Our numerical findings show that there is a trade-off between computational time and solution quality.
\end{abstract}

\textcolor{blue}{
\section*{Notation}
The main notation is listed below while other symbols are defined throughout the paper as needed. A subscript $t/\gamma/h/k$ in the notation refers to the corresponding values in the $t^{th}$ time stage/ $\gamma^{th}$ long-term scenario/ $h^{th}$ operating condition/ $k^{th}$ market scenario. Superscript/subscript $(\cdot)$ stands for the existing (E/{\it e}) and candidate (C/{\it c}) generation units, respectively. In addition, superscripts Conv and WP stand for conventional and wind power units, respectively.
\subsection{Sets and Indices}
\begin{ldescription}{$xxxxxxxx$}
\setlength{\itemsep}{2.5pt}
\item [$c \in \mathcal{C}$]  Set of candidate generation units.
\item [$d \in \mathcal{D}$]  Set of demands.
\item [$e \in \mathcal{E}$]  Set of existing generation units.
\item [$h \in \mathcal{H}$]  Set of wind-load operating conditions.
\item [$(k,k') \in \mathcal{K}$]  Set of short-term market scenarios.
\item [$r \in \mathcal{R}$]  Set of rival generation units.
\item [$(t,\tau)\in \mathcal{T}$] Set of time stages in the planning horizon.
\item [$(\gamma,\gamma') \in \mathcal{G}$]  Set of long-term scenarios.
\end{ldescription}
\subsection{Parameters}
\begin{ldescription}{$xxxxxxx$}
\setlength{\itemsep}{2.5pt}
\item [$a_{t}$] Amortization rate [\%].
\item [$b_{tkd}^{\text{D}}$] Utility of demand $d$ [\$/MWh]. 
\item [$c_{t\gamma (\cdot)}^{(\cdot)}$] Marginal cost of generation unit $(\cdot)$ [\$/MWh].
\item [$c_{t\gamma c}^{\text{inv}}$] Capital cost of candidate unit $c$ [\$/MW]. 
\item [$c_{t\gamma k r}^{\text{R}}$] Offering price of rival unit $r$ [\$/MWh]. 
\item [$\text{DF}_{t}$] Discount factor [\%].
\item [$\overline{I}_{t}$] Investment budget [\$]. 
\item [$K_{h(\cdot)}^{(\cdot),\text{CF}}$] Capacity factor of wind power unit $(\cdot)$ [p.u.].
\item [$K_{hd}^{\text{D}}$] Demand factor of demand $d$ [p.u.].
\item [$N_{h}^{\text{OC}}$] Weight of operating condition $h$ [h].
\item [$\overline{P}_{t\gamma h k r}^{\text{R}}$] Offering quantity of rival unit $r$ [MW]. 
\item [$\overline{P}_{t\gamma d}^{\text{D}}$] Maximum load of demand $d$ [MW]. 
\item [$X_{e}^{\text{E}}$] Installed capacity of existing unit $e$ [MW]. 
\item [$\overline{X}_{c}^{\text{C}}$] Maximum capacity of candidate unit $c$ [MW]. 
\item [$\pi_{\gamma}^{\text{LT}}/\pi_{k}^{\text{MS}}$] Probability of long-term/market scenario [-].
\item [$\chi^{\text{SoS}}$] Security of supply factor [p.u.].
\end{ldescription}
\subsection{Decision variables}
\begin{ldescription}{$xxxxxxx$}
\setlength{\itemsep}{2.5pt}
\item [$\overline{P}_{t\gamma h k (\cdot)}^{(\cdot)}$] Offering quantity of unit $(\cdot)$ [MWh].
\item [$P_{t\gamma h k (\cdot)}^{(\cdot)}$] Dispatch quantity of unit $(\cdot)$ [MWh].
\item [$P_{t\gamma hkd}^{\text{D}}$] Dispatch quantity of demand unit $d$ [MWh].
\item [$P_{t\gamma hk r}^{\text{R}}$] Dispatch quantity of rival unit $r$ [MWh].
\item [$X_{t\gamma c}^{\text{C}}$] Capacity of candidate unit $c$ [MW].
\item [$\beta_{t \gamma h k (\cdot)}^{(\cdot)}$] Offering price of unit $(\cdot)$ [\$/MWh].
\item [$\lambda_{t \gamma h k}$] Market-clearing price [\$/MWh]. 
\end{ldescription}
}

\section{Introduction}

\textcolor{blue}{
Among various decision-making problems in power systems, generation investment problems are one of the most complex to tackle from the computational point of view. They need to comprehensively account for different sources of uncertainty, including short-term (e.g., renewable production) and long-term (e.g., demand growth) \cite{cite_1}. They are even more complicated in a market environment due to uncertainty induced by market participation strategies of competing producers \cite{cite_2}. The computational burden of these problems is further increased for a price-maker\footnote{Unlike price-takers, a price-maker producer is capable of altering market equilibrium outcomes to its own benefit by making strategic offering decisions.} (strategic) producer since it requires a closed-loop system to model the impacts of its strategic decisions on market outcomes \cite{cite_5}-\cite{cite_4}.
}

\textcolor{blue}{
One natural approach to model this closed-loop system is to use bilevel programming \cite{cite_6}, which itself is a computationally demanding framework. There is an extensive literature exploring the use of the bilevel problems for the market-based generation investment -- see \cite{cite_2} for a thorough survey. The bilevel investment problems normally recast as mixed-integer linear programming (MILP) problems \cite{cite_7,cite_9}. Thus, they are prone to certain computational limitations and generally underperform in case of realistically sized power networks. There are two general practices to reduce computational complexity of bilevel investment problems: (\textit{i}) to introduce simplifying assumptions, e.g., ignore the dynamic (multi-year) representation of investment decisions or discard considering all short- and long-term uncertainty sources, and (\textit{ii}) to implement decomposition techniques. The second practice is generally more preferable since it allows for computing more informed investment solutions in polynomial time.  
}

\textcolor{blue}{
Decomposition methods for MILP problems generally fall into two categories: stage-based
methods, e.g., Benders decomposition \cite{c88} and its variations, and scenario-based, e.g., consensus alternating direction method of multipliers (consensus-ADMM) \cite{cite_13}, that is also referred to as progressive hedging \cite{cite_10}-\cite{cite_12}. The benefit of the former methods is that the optimally of the solution might be controlled over iterations through two bounds provided by a master problem and a set of sub-problems. However, the computational complexity of the master problem increases due to new cuts added at each iteration. The decomposition methods based on ADMM, instead, distribute the computational load among subproblems proportionally, and their complexity does not increase over iterations.  As a shortcoming, there is no guarantee that they necessarily converge to the global optimum in case of MILP problems. However, recent developments propose provable performance guarantees for such problems \cite{cite_14}-\cite{cite_16}.
}

\textcolor{blue}{
This paper proposes a scenario-based distributed algorithm based on consensus-ADMM to solve strategic investment problems with extensive representation of both long- and short-term uncertainties and multi-stage planning horizon. Unlike traditional algorithms in \cite{cite_11} and \cite{cite_12}, the proposed algorithm relaxes non-anticipativity conditions of both long- and short-term decision trees, thus splitting the original bilevel problem into a set of smaller bilevel problems with significantly lower computational needs. Using the framework of \cite{cite_16}, we prove the existence of the global bound on the optimal solution of the original bilevel problem. We then introduce an alternative local bound based on the tightness of nodes of short- and long-term decisions trees. The two bounds are to converge over iterations allowing for a practical performance guarantee: if the gap between the bounds closes at the last iteration, the algorithm provides the global optimal solution. 
}




The remainder of this paper is organized as follows. Section \ref{Sec_2} describes the considered strategic generation investment model and its reformulation as a MILP problem. Section \ref{Sec_3} explains the proposed consensus-ADMM algorithm and bounds on the optimal solution. Section \ref{Sec_4} illustrates the application of the algorithm and its ability to reach the global optimum. Section \ref{Sec_5} concludes the paper.

\section{Strategic Investment Problem} \label{Sec_2}

\subsection{Uncertainty and Decision Trees}
\tikzstyle{Node_LT} = [circle, minimum width=0.05cm, minimum height=0.05cm,text centered, draw=black]
\begin{figure}
\centering
\resizebox{8.8cm}{!}{%
\begin{tikzpicture}
\draw[thick] (0,0) -- node[midway, below, sloped] {$\dots$} (2,1);
\draw[thick] (0,0) -- node[midway, below, sloped] {$\dots$}  (2,-1);
\draw[thick] (2,1) -- node[midway, below, sloped] {$\dots$}  (4,1.5);
\draw[thick] (2,1) -- node[midway, below, sloped] {$\dots$}  (4,0.5);
\draw[thick] (2,-1) -- node[midway, below, sloped] {$\dots$}  (4,-0.5);
\draw[thick] (2,-1) -- node[midway, below, sloped] {$\dots$}  (4,-1.5);
\draw[thick] (5,2) -- node[midway, below, sloped] {$\dots$}  (7,2.5);
\draw[thick] (5,2) -- node[midway, below, sloped] {$\dots$}  (7,1.5);
\draw[thick] (5,0) -- node[midway, below, sloped] {$\dots$}  (7,0.5);
\draw[thick] (5,0) -- node[midway, below, sloped] {$\dots$}  (7,-0.5);
\draw[thick] (5,-2) -- node[midway, below, sloped] {$\dots$}  (7,-2.5);
\draw[thick] (5,-2) -- node[midway, below, sloped] {$\dots$}  (7,-1.5);
\node[draw=white ] at (1.5,0) {$\dots$};
\node[draw=white ] at (4.5,0) {$\dots$};
\node[draw=white ] at (5.5,1) {$\dots$};
\node[draw=white ] at (5.5,-1) {$\dots$};
\draw [dashed] (0,0) -- (0,-3);
\draw [dashed] (2,1) -- (2,-3);
\draw [dashed] (4,1.5) -- (4,-3);
\draw [dashed] (5,2) -- (5,-3);
\draw [dashed] (7,2.5) -- (7,-3);
\node[draw=white ] at (0,-3.9) {$(X_{t_{1}\gamma c}^{\text{C}})$};
\node[draw=white ] at (2,-3.9) {$(X_{t_{2}\gamma c}^{\text{C}})$};
\node[draw=white ] at (4,-3.9) {$(X_{t_{3}\gamma c}^{\text{C}})$};
\node[draw=white ] at (5.49,-3.9) {$(X_{t_{(N^{t}-1)}\gamma c}^{\text{C}})$};
\node[draw=white ] at (7.1,-3.9) {$(X_{t_{N^{t}}\gamma c}^{\text{C}})$};
\node[draw=white, xshift=0pt] at (0,-3.3) { \footnotesize Stage 1};
\node[draw=white, xshift=0pt] at (2,-3.3) { \footnotesize Stage 2};
\node[draw=white, xshift=0pt] at (4,-3.3) { \footnotesize Stage 3};
\node[draw=white, xshift=0pt] at (5.5,-3.3) { \footnotesize Stage $(N^{t}-1)$};
\node[draw=white, xshift=0pt] at (7.2,-3.3) { \footnotesize Stage $N^{t}$};
\node[draw=white ] at (7.5,2.5) {$\gamma_{1}$};
\node[draw=white ] at (7.5,1.5) {$\gamma_{2}$};
\node[draw=white ] at (7.6,-2.5) {$\gamma_{N^{\gamma}}$};
\node[draw=white ] at (7.9,-1.5) {$\gamma_{(N^{\gamma}-1)}$};
\node[rectangle, draw=black,fill=black!25,rounded corners=2.5] at (0.5,0.25) {};
\node[draw=white ] at (0.5,0.65) {$h_{1}$};
\node[draw=white ] at (1.5,1.15) {$h_{N^{h}}$};
\node[rectangle, draw=black,fill=black!25,rounded corners=2.5] at (1.5,0.75) {};
\node[rectangle, draw=black,fill=black!25,rounded corners=2.5] at (0.5,-0.25) {};
\node[rectangle, draw=black,fill=black!25,rounded corners=2.5] at (1.5,-0.75) {};
\node[rectangle, draw=black,fill=black!25,rounded corners=2.5] at (2.5,1.13) {};
\node[rectangle, draw=black,fill=black!25,rounded corners=2.5] at (3.5,1.37) {};
\node[rectangle, draw=black,fill=black!25,rounded corners=2.5] at (2.5,0.87) {};
\node[rectangle, draw=black,fill=black!25,rounded corners=2.5] at (3.5,0.63) {};
\node[rectangle, draw=black,fill=black!25,rounded corners=2.5] at (2.5,-0.87) {};
\node[rectangle, draw=black,fill=black!25,rounded corners=2.5] at (3.5,-0.63) {};
\node[rectangle, draw=black,fill=black!25,rounded corners=2.5] at (2.5,-1.13) {};
\node[rectangle, draw=black,fill=black!25,rounded corners=2.5] at (3.5,-1.37) {};
\node[rectangle, draw=black,fill=black!25,rounded corners=2.5] at (5.5,2.13) {};
\node[rectangle, draw=black,fill=black!25,rounded corners=2.5] at (6.5,2.37) {};
\node[rectangle, draw=black,fill=black!25,rounded corners=2.5] at (5.5,1.87) {};
\node[rectangle, draw=black,fill=black!25,rounded corners=2.5] at (6.5,1.63) {};
\node[rectangle, draw=black,fill=black!25,rounded corners=2.5] at (5.5,0.13) {};
\node[rectangle, draw=black,fill=black!25,rounded corners=2.5] at (6.5,0.37) {};
\node[rectangle, draw=black,fill=black!25,rounded corners=2.5] at (5.5,-0.13) {};
\node[rectangle, draw=black,fill=black!25,rounded corners=2.5] at (6.5,-0.37) {};
\node[rectangle, draw=black,fill=black!25,rounded corners=2.5] at (5.5,-2.13) {};
\node[rectangle, draw=black,fill=black!25,rounded corners=2.5] at (6.5,-2.37) {};
\node[rectangle, draw=black,fill=black!25,rounded corners=2.5] at (5.5,-1.87) {};
\node[rectangle, draw=black,fill=black!25,rounded corners=2.5] at (6.5,-1.63) {};
\node[circle, draw=black,fill=white] at (0,0) {};
\node[circle, draw=black,fill=white] at (2,1) {};
\node[circle, draw=black,fill=white] at (2,-1) {};
\node[circle, draw=black,fill=white] at (4,1.5) {};
\node[circle, draw=black,fill=white] at (4,0.5) {};
\node[circle, draw=black,fill=white] at (4,-0.5) {};
\node[circle, draw=black,fill=white] at (4,-1.5) {};
\node[circle, draw=black,fill=white] at (5,0) {};
\node[circle, draw=black,fill=white] at (7,0.5) {};
\node[circle, draw=black,fill=white] at (7,-0.5) {};
\node[circle, draw=black,fill=white] at (5,2) {};
\node[circle, draw=black,fill=white] at (7,2.5) {};
\node[circle, draw=black,fill=white] at (7,1.5) {};
\node[circle, draw=black,fill=white] at (5,-2) {};
\node[circle, draw=black,fill=white] at (7,-2.5) {};
\node[circle, draw=black,fill=white] at (7,-1.5) {};

\draw[->,>=stealth] (-0.5,-3) --  node[midway, below, sloped,yshift=-35] {Planning horizon}  (7.5,-3);

\draw[->,>=stealth] (12.85,-0.7) --  (13.3,-0.7);

\draw[thick,line width=0.3] (9.3,-0.7) -- (9.75,-0.7);

\node[draw=white] at (11.3,-0.7) {\text{Operating condition}};
\draw [thick, line width=0.3] (1.5,-4.5) -- (2.1,-4.5);
\draw [->,>=stealth, line width=0.3] (4.9,-4.5) -- (5.5,-4.5);
\draw [dashed] (6.5,2.5) -- (6.5,3);
\draw [dashed, line width=0.3] (8,2) -- (14.5,2);
\draw [dashed, line width=0.3] (8,-1) -- (14.5,-1);
\draw [dashed, line width=0.3] (8,2) -- (8,-1);
\draw [dashed, line width=0.3] (14.5,2) -- (14.5,-1);
\draw [dashed, line width=0.3] (9.2,0.85) -- (9.2,-0.05);
\draw [dashed, line width=0.3] (12.7,1.15) -- (12.7,-0.05);
\draw [dashed, line width=0.3] (6.5,3.0) -- (9.5,3.0);
\draw [dashed, line width=0.3] (9.5,3.0) -- (9.5,2.0);
\node[draw=white ] at (12.7,-0.3) {\text{Market clearing}};
\node[draw=white ] at (9.2,-0.3) {$\overline{\mathbf{P}}, \boldsymbol{\beta}$};
\draw[thick] (9.2,0.85) -- (12.7,1.15);
\draw[thick] (9.2,0.85) -- (12.7,0.55);
\node[draw=white ] at (12.2,0.85) {$\dots$};
\node[draw=white ] at (13.2,1.15) {$k_{1}$};
\node[draw=white ] at (13.34,0.55) {$k_{N^{k}}$};
\node[regular polygon, regular polygon sides=3, inner sep=0pt, minimum size=10pt, draw=black,fill=white] at (9.2,0.85) {};
\node[regular polygon, regular polygon sides=3, inner sep=0pt, minimum size=10pt, draw=black,fill=white] at (12.7,1.15) {};
\node[regular polygon, regular polygon sides=3, inner sep=0pt, minimum size=10pt, draw=black,fill=white] at (12.7,0.55) {};
\draw [->,>=stealth, line width=0.3] (8.7,-0.05) -- (13.7,-0.05);
\node[draw=white ] at (3.5,3.8) {\text{Long-term decision tree}};
\draw [decorate,decoration={brace,amplitude=10pt},line width=0.35mm, xshift=-4pt,yshift=0pt]
(-0.5,3.2) -- (7.5,3.2);
\node[draw=white ] at (11.2,3.8) {\text{Short-term decision tree}};
\draw [decorate,decoration={brace,amplitude=10pt},line width=0.35mm, xshift=-4pt,yshift=0pt]
(8.2,3.2) -- (14.7,3.2);

\node[circle, draw=black,fill=white] at (9.1,-2) {};
\node[draw=white ] at (11.7,-2) {- nodes for long-term scenarios};
\node[rectangle, draw=black,fill=black!25,rounded corners=2.5] at (9.1,-2.5) {};
\node[draw=white ] at (11.01,-2.5) {- operating conditions};
\node[regular polygon, regular polygon sides=3, inner sep=0pt, minimum size=10pt, draw=black,fill=white] at (9.1,-3) {};
\node[draw=white ] at (11.49,-2.96) {- nodes for market scenarios};
\node[draw=white ] at (10.62,-3.5) {$X_{t\gamma c}^{\text{C}}$ - investment decisions};
\node[draw=white ] at (11.21,-4.0) {$\overline{\mathbf{P}}, \boldsymbol{\beta}$ - market participation strategy};

\end{tikzpicture}
}
\caption{Long- and short-term decision trees of the strategic producer.} \label{fig:sequences}
\end{figure}
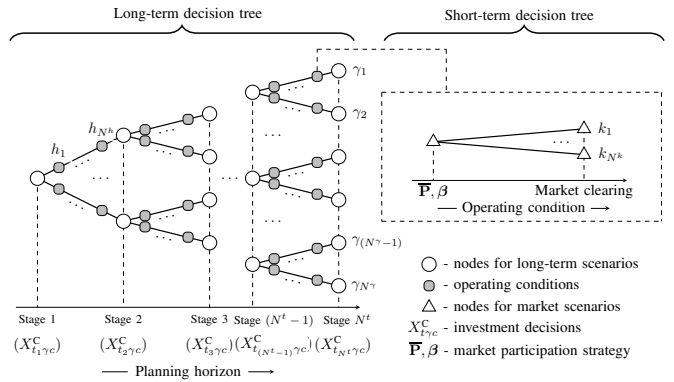

Investment decisions in power systems are subject to a wide range of uncertainties. To support informed decisions of the strategic producer, we account for the following uncertainty sources. First, the hourly variability of system load and wind power is considered through the finite set of operating conditions $\mathcal{H}$. Each condition is given by wind and load power factors and corresponding weight of that condition. The weight of each condition is the number of hours during the investment period represented by that condition. Second, short-term uncertainty is given by a set of market scenarios $\mathcal{K}$ describing the variability of price offering strategies of rival producers and demands.  Finally, long-term uncertainty set $\mathcal{G}$ contains the ambiguity of investment cost, demand growth and rivals' investment decisions. In this work, we rely on scenario representation of short- and long-term uncertainties.

Short- and long-term uncertainties shape the decision-making process of the producer as illustrated in Fig. \ref{fig:sequences}. At each time stage of the planning horizon, the producer decides investment $X_{t\gamma c}^{\text{C}}$ in candidate conventional and wind power units. Inside each investment period, it needs to decide participation strategy expressed through offering quantities $\overline{\mathbf{P}}$ and prices $\boldsymbol{\beta}$ for existing and candidate units.

\subsection{Bilevel Problem Formulation}
\tikzstyle{BOX} = [rectangle, rounded corners = 0, minimum width=180, minimum height=40,text centered, draw=black, fill=gray!2, line width=0.1mm]
\begin{figure}[]
\center
\begin{tikzpicture}

\node [align=center] (UL) [BOX] {\small{Upper-level problem  (\ref{UL_obj})-(\ref{price_limits})} \\ \small{(Expected profit maximization)}};
\node [align=center] (LL) [BOX, below of = UL, yshift=-50] {\small{Lower-level problems (\ref{LL_obj})-(\ref{candidate_cap_wp}) $ \forall t, \gamma, k, h$} \\ \small{(Market clearing)}};
\draw[transform canvas={xshift=-0.5cm},->,line width=0.2mm] (UL) -- node[midway, left, text width=2cm] {\scriptsize
$\overline{P}_{t\gamma h k c}^{\text{C}}$ , $\overline{P}_{t\gamma h k e}^{\text{E}}$, $\beta_{t \gamma h k e}^{\text{E}}$ , $\beta_{t \gamma h k c}^{\text{C}}$
}  (LL) ;
\draw[transform canvas={xshift=0.5cm},->,line width=0.2mm] (LL) -- node[midway, right, text width=2cm] {
\scriptsize
$P_{t\gamma h k c}^{\text{C}}$ , $P_{t\gamma h k e}^{\text{E}}$, $\lambda_{t \gamma h k}$
}  (UL) ;
\end{tikzpicture}
\caption{Bilevel structure of the investment problem and the interactions
between the upper- and lower-level problems.}
\label{bilevel_scheme}
\vspace{-0.5cm}
\end{figure}
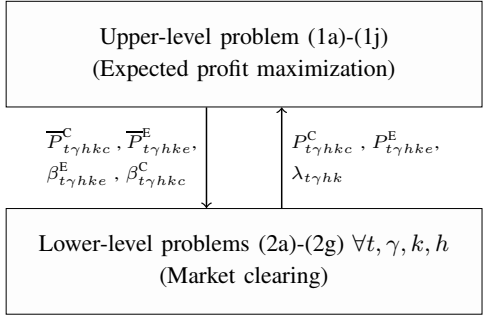
The proposed bilevel problem consists of an upper-level (UL) problem and a set of lower-level (LL) problems as depicted in Fig. \ref{bilevel_scheme}. The UL problem maximizes the expected profit of strategic power producer throughout the planning horizon by computing optimal investment and market participation decisions. Using this bilevel setup, the strategic producer anticipates the market clearing outcomes in the LL problems as a function of its strategic decisions made in the UL problem. The LL problems are specified for each time stage of the planning horizon, short-term scenario, long-term scenario and operating condition. These problems are interconnected in the sense that LL optimization problems are treated as constraints for the UL problem. The investment and participation decisions made in the UL problem affect the outcome of the LL problems that provide market prices and dispatch quantities that, in turn, affect the expected profit in the UL problem. The UL problem writes as the following multi-stage stochastic problem:
\begin{subequations} \label{UL_problem}
\begin{align}
  &\text{max}_{\Delta^{\text{UL}}} \  \sum_{t \in \mathcal{T}} \text{DF}_{t} \Big\{ \sum_{\gamma \in \mathcal{G}} \pi_{\gamma}^{\text{LT}} \Big[\sum_{h \in \mathcal{H}} N_{h}^{\text{OC}} \Big<\sum_{k \in \mathcal{K}}
  \pi_{k}^{\text{MS}} \nonumber\\
  &\Big(
  \sum_{e \in \mathcal{E}} (\lambda_{t \gamma h k} - c_{t \gamma e}^{\text{E}}) P_{t \gamma h k e}^{\text{E}} +
  \sum_{c \in \mathcal{C}} (\lambda_{t \gamma h k} - c_{t \gamma c}^{\text{C}}) P_{t \gamma h k c}^{\text{C}}
  \Big)\Big> \nonumber\\
  & - a_{t} \sum_{c \in \mathcal{C}} c_{t \gamma c}^{\text{inv}} \sum_{\tau \in \mathcal{T} \atop \tau \leq t} X_{\tau \gamma  c}^{\text{C}} \Big]\Big\} \label{UL_obj} \\
  &\text{s.t.} \quad X_{t\gamma c}^{\text{C}} = X_{t\gamma' c}^{\text{C}} \quad \forall t,(\gamma,\gamma')\in \overline{\mathbf{G}}_{t},c, \label{non_anticip_cond} \\
  &0 \leq X_{t\gamma c}^{\text{C}} \leq \overline{X}_{c}^{\text{C}} \quad \forall t,\gamma,c,\label{install_cap}\\
  & \sum_{c \in \mathcal{C}} c_{t \gamma c}^{\text{inv}}  X_{t\gamma c}^{\text{C}} \leq \overline{I}_{t} \quad \forall t,\gamma, \label{budget} \\
  & \sum_{c \in \mathcal{C}} \overline{P}_{t\gamma h k c}^{\text{C}} + \sum_{e \in \mathcal{E}} \overline{P}_{t\gamma h k e}^{\text{E}} + \sum_{r \in \mathcal{R}} \overline{P}_{t\gamma h k r}^{\text{R}} \geq \nonumber\\
  &\quad\quad\quad\quad\quad\quad\quad\quad\chi^{\text{SoS}}  \sum_{d \in \mathcal{D}} \overline{P}_{t\gamma d}^{\text{D}} K_{h}^{\text{DF}} \quad  \forall t, \gamma, h, k, \label{Security_of_supply} \\
  & 0 \leq \overline{P}_{t\gamma h k c}^{\text{C}} \leq X_{t\gamma c}^{\text{C}} K_{hc}^{\text{CF}} \quad \forall t,\gamma,h,k,c \in \mathcal{C}^{\text{WP}}, \label{offer_limits_first} \\
  &0 \leq \overline{P}_{t\gamma h k c}^{\text{C}} \leq X_{t\gamma c}^{\text{C}} \quad \forall t,\gamma,h,k,c \in \mathcal{C}^{\text{Conv}}, \\
  &0 \leq \overline{P}_{t\gamma h k e}^{\text{E}} \leq X_{t\gamma e}^{\text{E}} K_{he}^{\text{CF}} \quad \forall t,\gamma,h,k, e \in \mathcal{E}^{\text{WP}},\\
  &0 \leq \overline{P}_{t\gamma h k e}^{\text{E}} \leq X_{t\gamma e}^{\text{E}} \quad \forall t,\gamma,h,k,e \in \mathcal{E}^{\text{Conv}}, \label{offer_limits_last}\\
  & \beta_{t \gamma h k e}^{\text{E}}, \beta_{t \gamma h k c}^{\text{C}} \geq 0 \quad \forall t,\gamma,h,k,c,e, \label{price_limits}
\end{align}
\end{subequations}
where $\Delta^{\text{UL}} \in \{X_{t\gamma c}^{\text{C}}, \overline{P}_{t\gamma h k c}^{\text{C}}, \overline{P}_{t\gamma h k e}^{\text{E}}, \beta_{t \gamma h k e}^{\text{E}}, \beta_{t \gamma h k c}^{\text{C}}\}$ is the set of strategic producer's decision variables, comprising investment decisions in candidate units and offering quantities and prices for both existing and candidate units. The UL objective function (\ref{UL_obj}) is discounted expected profit from operations of existing and candidate units subtracting investment costs. Constraints (\ref{non_anticip_cond}) are non-anticipativity conditions on investment decisions enforced at each time stage by incidence matrix $\overline{\mathbf{G}}_{t}$. The matrix $\overline{\mathbf{G}}_{t}$ ensures that investment decisions in adjacent scenarios $\gamma$ and $\gamma'$ at time stage $t$ are identical for all possible realizations of long-term uncertainty set $\mathcal{G}$ at time stages following $t$. Inequalities (\ref{install_cap}) and (\ref{budget}) limit the installed capacity of candidate units and associated expenses with upper bounds. Regulatory constraints (\ref{Security_of_supply}) enforcing security of supply prevent the strategic producer from causing capacity shortage in the system. Finally, a set of constraints (\ref{offer_limits_first})-(\ref{price_limits}) defines bounds on the supply functions, i.e., on offering power quantities and associated prices for each existing and candidate generation unit. The dispatch quantities and market clearing prices are treated as parameters in the UL problem that are obtained by solving the following set of the LL market clearing problems:
\begin{subequations} \label{LL_problem}
\begin{align}
&\Big\{\text{max}_{\Delta^{\text{LL}}} \quad \sum_{d \in \mathcal{D}} b_{tkd}^{\text{D}} P_{t\gamma hkd}^{\text{D}} - \sum_{r \in \mathcal{R}} c_{t\gamma kr}^{\text{R}} P_{t\gamma hk r}^{\text{R}} - \nonumber \\
&\quad\quad\quad\quad\quad\quad\sum_{c \in \mathcal{C}} \beta_{t \gamma h k c}^{\text{C}} P_{t\gamma h k c}^{\text{C}} -  \sum_{e \in \mathcal{E}} \beta_{t \gamma h k e}^{\text{E}} P_{t\gamma h k e}^{\text{E}} \label{LL_obj} \\
&\text{s.t.} \quad \sum_{r \in \mathcal{R}} P_{t\gamma hk r}^{\text{R}} + \sum_{c \in \mathcal{C}} P_{t\gamma hk c}^{\text{C}} + \sum_{e \in \mathcal{E}} P_{t\gamma hk e}^{\text{E}} \nonumber\\
&\quad\quad\quad\quad\quad\quad\quad\quad\quad\quad\quad- \sum_{d \in \mathcal{D}} P_{t\gamma hkd}^{\text{D}} = 0 : \lambda_{t \gamma h k} \label{DA_PB} \\
&0 \leq P_{t\gamma hkd}^{\text{D}} \leq \overline{P}_{t\gamma d}^{\text{D}} K_{hd}^{\text{DF}} \quad \forall d \label{demand_cap}\\
&0 \leq P_{t\gamma hk r}^{\text{R}} \leq \overline{P}_{t\gamma h k r}^{\text{R}} \quad \forall r \in \mathcal{R}^{\text{Conv}} \label{rival_cap}\\
&0 \leq P_{t\gamma hk r}^{\text{R}} \leq \overline{P}_{t\gamma h k r}^{\text{R}} K_{hr}^{\text{R,CF}} \quad \forall r \in \mathcal{R}^{\text{WP}} \label{rival_cap_wp}\\
&0 \leq P_{t\gamma h k e}^{\text{E}} \leq \overline{P}_{t\gamma h k e}^{\text{E}} \quad \forall e \\
&0 \leq P_{t\gamma h k c}^{\text{C}} \leq \overline{P}_{t\gamma h k c}^{\text{C}} \quad \forall c \quad\Big\} \forall t,\gamma,h,k, \label{candidate_cap_wp}
\end{align}
\end{subequations}
where $\Delta^{\text{LL}} \in \{P_{t\gamma hkd}^{\text{D}}, P_{t\gamma hk r}^{\text{R}}, P_{t\gamma h k c}^{\text{C}}, P_{t\gamma h k e}^{\text{E}}\}$ is the set of primal LL decision variables that includes dispatch of generation and load units for each time stage, long-term scenario, short-term scenario and operating condition. In addition, $\lambda_{t \gamma h k}$ is a market clearing price that is obtained as dual variable of (\ref{DA_PB}). The LL objective function (\ref{LL_obj}) represents the market social welfare, subject to power balance (\ref{DA_PB}) and dispatch limits of generation and load units (\ref{demand_cap})-(\ref{candidate_cap_wp}). All UL variables are treated as parameters within the LL problems, which makes the LL problems linear and convex.

\textcolor{blue}{To derive a single level equivalent, the lower-level problems \eqref{LL_problem} are replaced with their Karush-–Kuhn–-Tucker (KKT) conditions. The non-linear terms in (\ref{UL_obj}), i.e., product of dual prices and dispatch quantities of existing and candidate units, are replaced with their exact linear equivalents as explained in \cite{cite_2}. Besides, the linear equivalents of the complementarity slackness conditions are obtained using special ordered set of type 1 (SOS1) variables as explained in \cite{cite_17}. As a result, the bilevel problem is recast as a single-level MILP problem.}


\section{Proposed Consensus-ADMM Algorithm} \label{Sec_3}
\subsection{Algorithm Description}
Decomposing the single-level equivalent of strategic investment problem (\ref{UL_problem})-(\ref{LL_problem}) per long-term scenario by relaxing non-anticipativity constraints (\ref{non_anticip_cond}) would result in a number of sub-problems corresponding to the size of scenario set $\mathcal{G}$. Resulting sub-problems are still stochastic problems due to short-term uncertainty accounted for in set $\mathcal{K}$, and themselves might be still difficult to solve. Thus, the relaxation of the long-term decision tree might not be sufficient to reduce computational complexity of the problem. Our algorithm suggests to relax decision trees associated with \textit{both} long-term and short-term uncertainties, such that the resulting sub-problems become deterministic, requiring less computational effort to solve.

Let $\mathbf{X}_{t \gamma k}$ be a vector of investment decisions in a set of investment options $\mathcal{C}$ at time stage $t$ that is specific for a pair of short- and long-term scenarios. Then, for particular long- and short-term scenarios $\gamma'$ and $k'$, the non-anticipativity constraints (\ref{non_anticip_cond}) are reformulated as follows:
\begin{align} \label{non_ant_alg}
         &\mathbf{X}_{t \gamma k} - \mathbf{\overline{X}}_{t \gamma' k'} = 0
         \quad \forall t,(k,k') \in \overline{\mathbf{K}}, (\gamma,\gamma') \in \overline{\mathbf{G}}_{t},
\end{align}
where $\mathbf{\overline{X}}_{t \gamma' k'}$ is a global variable which requires scenario-specific investment decisions to coincide according to the conditions enforced by long- and short-term non-anticipativity matrices $\overline{\mathbf{G}}_{t}$ and $\overline{\mathbf{K}}$. Unlike $\overline{\mathbf{G}}_{t}$, matrix $\overline{\mathbf{K}}$ states that short-term scenario-specific investment solutions have to coincide at all time stages of the planning horizon. By relaxing (\ref{non_ant_alg}), the amount of sub-problems is now defined by a number of long- and short-term scenarios. Let us then denote a coefficient vector and a vector of all decision variables of each sub-problem by $\mathbf{c}_{t\gamma k}$ and $\mathbf{x}_{t}$, respectively. Vector $\mathbf{x}_{t} \in \mathcal{Q}_{t \gamma k},$ where $\mathcal{Q}_{t \gamma k}$ is a time- and scenario-specific non-convex feasible set of each sub-problem. By $\mu_{t\gamma k}$ we denote a dual variable of (\ref{non_ant_alg}). Then, the proposed iterative algorithm writes as follows:
\begin{align}
\mathbf{X}_{t\gamma k}^{\nu} \leftarrow\quad& \underset{\mathbf{x}_{t} \in \mathcal{Q}_{t \gamma k}}{\text{argmax}}
\Big\{
\sum_{t\in T}
\big(
\mathbf{c}_{t\gamma k}^{\top} \mathbf{x}_{t} -
\mu_{t\gamma k}^{\nu-1 \top} \mathbf{x}_{t} \nonumber\\
&\quad\quad\quad-\frac{\rho}{2} \norm{\mathbf{x}_{t} - \mathbf{X}_{t\gamma k}^{\nu-1}}_{2}^{2}
\big)
\Big\}, &&\forall \gamma,k, \label{step_first}\\
\mathbf{\overline{X}}_{t \gamma k}^{\nu} \leftarrow\quad&
\frac{\sum_{\gamma' \in \overline{\mathbf{G}}_{t} \atop k' \in \overline{\mathbf{K}}} \pi_{\gamma'}^{\text{LT}} \pi_{k'}^{\text{MS}} \mathbf{X}_{t\gamma' k'}^{\nu}}{\sum_{\gamma' \in \overline{\mathbf{G}}_{t} \atop k' \in \overline{\mathbf{K}}} \pi_{\gamma'}^{\text{LT}} \pi_{k'}^{\text{MS}}}, &&\forall t,\gamma,k, \label{step_second}\\
\mu_{t\gamma k}^{\nu} \leftarrow\quad& \mu_{t\gamma k}^{\nu-1} + \rho
\big(
\mathbf{X}_{t\gamma k}^{\nu} - \mathbf{\overline{X}}_{t \gamma k}^{\nu}
\big), &&\forall t,\gamma,k, \label{step_three}
\end{align}
where $\nu$ is an index of iterations. As the first step, investment decisions are obtained in sub-problems (\ref{step_first}) for each pair of long- and short-term scenario and a time stage. The objective function of (\ref{step_first}) is represented by the scenario-specific objective function (\ref{UL_obj}), augmented by two penalization terms. The first term results from augmenting (\ref{non_ant_alg}) into objective function of each sub-problem and aims at adjusting the investment solutions towards the mean of adjacent nodes, while the second proximal term drives the algorithm towards convergence. As the second step, the algorithm updates the global variable in (\ref{step_second}) as a probability-weighted average solution over adjacent scenarios, defined by matrices $\overline{\mathbf{G}}_{t}$ and $\overline{\mathbf{K}}$. Last step is a dual update according to (\ref{step_three}), where factor $\rho$ penalizes the deviation of specific investment decisions from the corresponding average solution. Convergence of the algorithm is verified over iterations by $\mathbf{g}_{t \gamma k}^{\nu}$ which indicates weather scenario-specific investment decisions coincide with respective global variable, such that:
\begin{align}
         & \mathbf{g}_{t \gamma k}^{\nu} \leftarrow |\mathbf{X}_{t\gamma k}^{\nu} - \mathbf{\overline{X}}_{t \gamma k}^{\nu}|, \quad \forall t, \gamma, k. \label{convergence_g}
\end{align}
Convergence is reached when \eqref{convergence_g} remains below a predefined tolerance $\epsilon$.

\subsection{Bounds and \textcolor{blue}{Performance Guarantee}}
\textcolor{blue}{Here we aim at introducing two bounds on the optimal value of objective function (\ref{UL_obj}) that provide a practical performance guarantee for the proposed algorithm. Leveraging the framework in \cite{cite_16},} we introduce the global upper bound GUB for multi-stage investment problem with relaxation of long- and short-term decision trees as follows. We start by introducing the following proposition.


\textbf{Proposition 1.}
{\it
By denoting a vector of optimal investment decisions as $\mathbf{\dot{X}}_{t\gamma k}^{\nu}$, the following condition holds for each iteration of the algorithm:
}
\begin{center}
$\sum_{\gamma \in \overline{\mathbf{G}}_{t} \atop k \in \overline{\mathbf{K}}} \pi_{\gamma}^{\text{LT}} \pi_{k}^{\text{MS}} \mu_{t\gamma k}^{\nu \top} \mathbf{\dot{X}}_{t\gamma k}^{\nu} = 0 \quad \forall t \in \mathcal{T}.$
\end{center}
\begin{proof}
It is provable by induction. Let consider iteration zero, in which dual update (\ref{step_three}) is defined as
\begin{center}
$\mu_{t\gamma k}^{0} = \rho
\big(
\mathbf{X}_{t\gamma k}^{0} - \mathbf{\overline{X}}_{t \gamma k}^{0}
\big), \ \forall t,\gamma,k.$
\end{center}

Then, by definition of $\mathbf{\overline{X}}_{t \gamma k}$, in expectation it rewrites as

\begin{align*}
&\sum_{\gamma \in \overline{\mathbf{G}}_{t} \atop k \in \overline{\mathbf{K}}} \pi_{\gamma}^{\text{LT}} \pi_{k}^{\text{MS}} \mu_{t\gamma k}^{0} = \rho \sum_{\gamma \in \overline{\mathbf{G}}_{t} \atop k \in \overline{\mathbf{K}}} \pi_{\gamma}^{\text{LT}} \pi_{k}^{\text{MS}}
\big(
\mathbf{X}_{t\gamma k}^{0} - \mathbf{\overline{X}}_{t \gamma k}^{0}
\big) = \\
& \rho \sum_{\gamma \in \overline{\mathbf{G}}_{t} \atop k \in \overline{\mathbf{K}}} \pi_{\gamma}^{\text{LT}} \pi_{k}^{\text{MS}}
\frac
{\sum_{\gamma' \in \overline{\mathbf{G}}_{t} \atop k' \in \overline{\mathbf{K}}} \pi_{\gamma'}^{\text{LT}} \pi_{k'}^{\text{MS}} (\mathbf{X}_{t\gamma k}^{0} - \mathbf{X}_{t \gamma' k}^{0})}
{\sum_{\gamma' \in \overline{\mathbf{G}}_{t} \atop k' \in \overline{\mathbf{K}}} \pi_{\gamma'}^{\text{LT}} \pi_{k'}^{\text{MS}} } = 0.
\end{align*}
By induction, the same holds for subsequent iterations.
\end{proof}
We now define  $D_{\gamma k}^{\nu}$ as an optimal solution to the following problem:
\begin{center}
$D_{\gamma k}^{\nu} = \underset{\mathbf{x}_{t} \in \mathcal{Q}_{t \gamma k}}{\text{max}} \sum_{t \in \mathcal{T}} \Big(
\mathbf{c}_{t\gamma k}^{\top} \mathbf{x}_{t} -
\mu_{t\gamma k}^{\nu \top} \mathbf{x}_{t}
\Big).$
\end{center}
Then, the global upper bound GUB is introduced with the following theorem.

\textbf{Theorem 1.}
{\it
By denoting the global optimal solution of the stochastic problem (\ref{UL_problem})-(\ref{LL_problem}) as $\dot{z}$, the following condition holds at each iteration of the algorithm:
}
\begin{center}
$\text{GUB} = \sum_{\gamma \in \overline{\mathbf{G}}_{t} \atop k \in \overline{\mathbf{K}}} \pi_{\gamma}^{\text{LT}} \pi_{k}^{\text{MS}} D_{\gamma k}^{\nu} \geq \dot{z}.$
\end{center}
\begin{proof}
From the definition of $D_{\gamma k}$,
\begin{center}
$D_{\gamma k} \geq \sum_{t \in \mathcal{T}} \Big(\mathbf{c}_{t\gamma k}^{\top} \mathbf{\dot{x}}_{t} - \mu_{t\gamma k}^{\nu \top} \mathbf{\dot{x}}_{t} \Big) \ \forall \gamma,k.$
\end{center}
Taking into account Proposition 1,
\begin{align*}
\sum_{\gamma \in \overline{\mathbf{G}}_{t} \atop k \in \overline{\mathbf{K}}} \pi_{\gamma}^{\text{LT}} \pi_{k}^{\text{MS}} D_{\gamma k}^{\nu} \geq &
\sum_{\gamma \in \overline{\mathbf{G}}_{t} \atop k \in \overline{\mathbf{K}}} \pi_{\gamma}^{\text{LT}} \pi_{k}^{\text{MS}} \sum_{t \in \mathcal{T}} \Big(\mathbf{c}_{t\gamma k}^{\top} \mathbf{\dot{x}}_{t} - \mu_{t\gamma k}^{\nu \top} \mathbf{\dot{x}}_{t} \Big)\\
\geq & \sum_{\gamma \in \overline{\mathbf{G}}_{t} \atop k \in \overline{\mathbf{K}}} \pi_{\gamma}^{\text{LT}} \pi_{k}^{\text{MS}} \sum_{t \in \mathcal{T}} \mathbf{c}_{t\gamma k}^{\top} \mathbf{\dot{x}}_{t} = \dot{z}.
\end{align*}
Thus, at any iteration $\nu$, $\text{GUB}^{\nu} \geq \dot{z}$.
\end{proof}

We then introduce a local upper bound denoted by UB, that is defined based on the tightness of the adjacent nodes of the relaxed short- and long-term decision trees. By fixing investment decisions to the ones provided by (\ref{step_first}), at each iteration the UB is computed as follows:
\begin{align*}
         & \text{UB}^{\nu} = \sum_{\gamma \in \overline{\mathbf{G}}_{t} \atop k \in \overline{\mathbf{K}}} \pi_{\gamma}^{\text{LT}} \pi_{k}^{\text{MS}} \Big[\underset{\mathbf{x}_{t} \in \mathcal{Q}_{t \gamma k} \atop \mathbf{x}_{t} \in \mathbf{X}_{t\gamma k}^{\nu}}{\text{max}} \sum_{t \in \mathcal{T}} \mathbf{c}_{t\gamma k}^{\top} \mathbf{x}_{t}\Big].
\end{align*}

\textcolor{blue}{The two bounds tend to a common basis since the nodes of short- and long-term trees get tighter over iterations. By definition, $\lim_{\mathbf{g}_{t \gamma k}\rightarrow 0}\text{UB} = \dot{z}$. As shown in \cite{cite_16}, with a proper tuning of penalty factor $\rho$, $\text{GUB} \rightarrow \dot{z}$. Consequently, we introduce a practical performance guarantee based on the gap between the two bounds. If $||\text{GUB-UB}||^2 = 0$ at the last iteration, the algorithm provides the global optimal solution of (\ref{UL_problem})-(\ref{LL_problem}), and this norm is nearly zero close to the optimum.
}


\vspace{-0.5cm}
\section{Simulation Results} \label{Sec_4}
We consider the instance of a moderate-scale power system to derive the optimal solution provided by the original extensive form of the stochastic MILP problem (\ref{UL_problem})-(\ref{LL_problem}). By extensive formulation solution, we mean the direct solution of (\ref{UL_problem})-(\ref{LL_problem}) without using decomposition. This optimal solution  is used as a benchmark for the proposed algorithm. Detailed data description and codes for all simulations are available in the online appendix of the paper \cite{cite_18}. The simulations are performed using CPLEX 12.1 under GAMS on an Intel Xeon processor E5-2680 with 8 cores clocking at 2.8 GHz and 128 GB of RAM.

The system initially consists of seven conventional generation units, five of which are rival units and two belong to the strategic producer. The total installed capacity of all generation units is 1500 MW. The load is represented by a single demand block of 1050 MW. The investment horizon consists of two time stages with three years in between. Three candidate technologies are available for investments: CCGT, coal, and wind power units, with investment costs increasing in that order. Investment budget is such that it is never binding in any scenario. The uncertainty of wind power production is described by five operating conditions, while demand factor is fixed to 1 across all operating conditions. The long-term uncertainty is characterized by three equiprobable demand growth scenarios, in which the demand at the second time period is 20\% higher, the same, or 20\% lower than that in the first period. Similarly, three market scenarios with equal probabilities are considered, such that the rival offering prices are 10\% higher, the same, or 10\% lower than the initial marginal costs. Further description of the test case is available in \cite{cite_18}.

We first directly solve the extensive MILP formulation of (\ref{UL_problem})-(\ref{LL_problem}) to obtain the optimal solution. Then, we apply the proposed consensus-ADMM algorithm with the relaxation of both long- and short-term decision trees. With this decomposition, the original problem is decomposed into a number of sub-problems, one per each pair of long- and short-term scenarios, i.e., nine sub-problems in this study. Table \ref{complexity} summarizes the complexity of two solution alternatives. Compared to the extensive formulation, the number of integer variables in each sub-problem is reduced by 89\%, so that their execution would require much less computational efforts.

\begin{table}[]
\centering
\caption{Computational complexity}
\label{complexity}
\begin{tabular}{lcc}
\Xhline{2\arrayrulewidth}
\multicolumn{1}{c}{Problem} & \multicolumn{1}{c}{\begin{tabular}[c]{@{}c@{}}Extensive\\ formulation\end{tabular}} & \multicolumn{1}{c}{\begin{tabular}[c]{@{}c@{}}Each ADMM\\ sub-problem\end{tabular}} \\
\Xhline{2\arrayrulewidth}
Number of variables & 9 648 & 1 180 \\
- Continues & 5 688 & 740 \\
- Integer & 3 960 & 440 \\
Number of constraints & 6 613 & 739 \\
\Xhline{2\arrayrulewidth}
\end{tabular}
\end{table}

\begin{table}[]
\centering
\caption{First-stage investment decisions [MW]}
\label{inv}
\tabcolsep=0.13cm
\begin{tabular}{M{1.3cm}M{1.3cm}M{1cm}M{1cm}M{1cm}N}
\Xhline{2\arrayrulewidth}
\multirow{2}{*}{Problem} & \multirow{2}{*}{\begin{tabular}[c]{@{}c@{}}Extensive\\ formulation\end{tabular}} & \multicolumn{3}{c}{ADMM} &\\[5pt]
\cline{3-6}
 &  & $\rho=10^2$ & $\rho=10^3$ & $\rho=10^5$ &\\[5pt]
\Xhline{2\arrayrulewidth}
CCGT & 0.0  & 0.0  & 0.0    & 0.2  &\\
Coal & 14.8 & 14.7 & 14.4 & 14.2 &\\
Wind & 88.2 & 88.3 & 88.5 & 88.6 &\\
\Xhline{2\arrayrulewidth}
\end{tabular}
\vspace{-0.5cm}
\end{table}

\begin{figure}[]
\center
\begin{tikzpicture}[thick,scale=1]
\pgfplotsset{compat=1.11, ymin=184, ymax=186, xmin=16, xmax=34, try min ticks=4}
\begin{axis}[
y tick label style={
/pgf/number format/.cd,
fixed,
fixed zerofill,
precision=1,
/tikz/.cd},
legend cell align=right,
legend columns=3,
legend style={fill=none,draw=none, at={(0.35,0.91)},
anchor=north,
font=\scriptsize},
name = ConvD_leg_pos,
font=\footnotesize,
width=0.5\textwidth,
height=0.2\textwidth,
xtick={5,10,15,20,25,30},
xticklabels={$\rho=10^2$,$\rho=10^3$,$\rho=10^5$,$\rho=10^2$,$\rho=10^3$,$\rho=10^5$},
ylabel={Expected profit [M$\$$]},
]
\addplot [color_ConvD, line width=0.1mm, smooth] table [x index = 0, y index = 1] {Direct_opt.dat};
\addlegendentry{Optimal}

\addplot[color=color_StochD,mark=none, line width=0.35mm] coordinates {
                  (19,0)
		(19,184.6)
                  (21,184.6)
                  (21,0)

	};
\addlegendentry{GUB}
\addplot[color={color_bilevel},mark=none, line width=0.2mm] coordinates {
                  (19,0)
         	(19,184.6)
                  (21,184.6)
                  (21,0)

         };
\addlegendentry{UB}
\addplot[color=color_StochD,mark=none, line width=0.35mm] coordinates {
                  (24,0)
                  (24,184.6)
                  (26,184.6)
                  (26,0)
};
\addplot[color={color_bilevel},mark=none, line width=0.2mm] coordinates {
                  (24,0)
                  (24,184.55)
                  (26,184.55)
                  (26,0)
};
\addplot[color=color_StochD,mark=none, line width=0.35mm] coordinates {
                  (29,0)
                  (29,185.4)
                  (31,185.4)
                  (31,0)
};
\addplot[color={color_bilevel},mark=none, line width=0.2mm] coordinates {
                  (29,0)
                  (29,184.5)
                  (31,184.5)
                  (31,0)
};
\end{axis}
\node at (3.6,-0.7) {\scriptsize{Penalty factor}};
\end{tikzpicture}
\caption{Bounds on the expected profit obtained for different values of penalty factor $\rho$}
\label{fig_gap}
\end{figure}
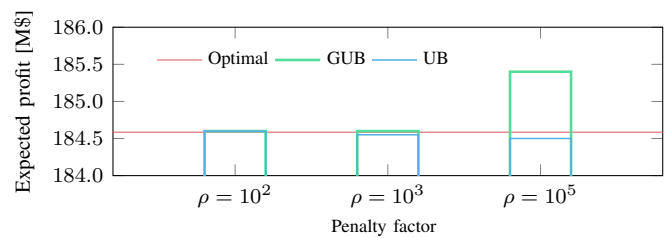

\begin{table}[]
\centering
\caption{Computational performance}
\label{sim_time}
\tabcolsep=0.13cm
\begin{tabular}{M{1.3cm}M{1.3cm}M{1cm}M{1cm}M{1cm}N}
\Xhline{2\arrayrulewidth}
\multirow{2}{*}{Problem} & \multirow{2}{*}{\begin{tabular}[c]{@{}c@{}}Extensive\\ formulation\end{tabular}} & \multicolumn{3}{c}{ADMM} &\\[3pt]
\cline{3-6}
 &  & $\rho=10^2$ & $\rho=10^3$ & $\rho=10^5$ &\\[3pt]
\Xhline{2\arrayrulewidth}
\begin{tabular}[c]{@{}c@{}}Number of \\ iterations\end{tabular} & - & 331 & 33 & 3 &\\[3pt]
Time {[}s{]} & 3624 & 1632 & 134 & 9 &\\[3pt]
\Xhline{2\arrayrulewidth}
\end{tabular}
\end{table}

In the optimal solution of extensive formulation, the expected profit amounts to $\$184.6$ million, while 14.8-MW of coal and 88.2-MW of wind power generation are built at the first time stage. The application of the proposed ADMM algorithm results in nearly the same investment solutions which depends on the setting of algorithm's parameters, as illustrated in Table \ref{inv}. With small values of penalty factor $\rho$, the solution in nearly identical to the optimal one with the slight difference explained by algorithm tolerance $\epsilon$, which is set to 0.5 MW. By increasing $\rho$, the solution deviates from the optimum in a sense that investment in the coal generation slightly decreases in favor of increased investment in stochastic wind generation. The ADMM algorithm estimates the expected profit in terms of bounds on the optimal solution as depicted in Fig. \ref{fig_gap}. It shows that the accuracy of the profit estimate reduces in $\rho$: for small $\rho$, both upper bounds coincide in the optimum, while with higher $\rho$ the estimate is distorted due to the increased gap between two bounds.

The computational performance of the proposed ADMM algorithm is compared with that of the extensive formulation in Table \ref{sim_time}. Among three values of $\rho$ tested, the simulation time for the ADMM algorithm is at most half as much as time required for the non-decomposed implementation, and it depends on the choice of penalty factor $\rho$. Small penalty factors result in more precise investment solutions but require more computational resources. Higher values of $\rho$, instead, drastically reduce the execution time, e.g., nine seconds against nearly an hour, at the expense of slight deviation from the global optimum. This way, by tuning the algorithm settings, a decision-maker can choose a trade-off between the quality of the solution and corresponding simulation time.

\begin{figure}[]
\center
\resizebox{7.5cm}{!}{%
\begin{tikzpicture}[thick,scale=1]
\pgfplotsset{compat=1.11, ymin=184, ymax=185.7, xmin=0, xmax=33, try min ticks=4}
\begin{axis}[
         y tick label style={
             /pgf/number format/.cd,
                 fixed,
                 fixed zerofill,
                 precision=1,
             /tikz/.cd},
xlabel={Number of iterations},
ylabel={Expected profit [M$\$$]},
ylabel near ticks,
legend cell align=right,
label style={font=\scriptsize},
tick label style={font=\scriptsize},
legend style={fill=none,draw=none, at={(0.28,1)},
anchor=north,
font=\tiny},
name = ConvD_leg_pos,
font=\scriptsize,
width=0.5\textwidth,
height=0.28\textwidth
]
\addplot [green, dashed, line width=0.35mm, smooth] table [x index = 0, y index = 1] {GUB_PH_LT_MS_10_5.dat};
\addlegendentry{\text{GUB}, $\rho=10^{5}$}
\addplot [green, line width=0.35mm, smooth] table [x index = 0, y index = 1] {UB_PH_LT_MS_10_5.dat};
\addlegendentry{\text{UB}, $\rho=10^{5}$}
\end{axis}
\begin{axis}[
hide axis,
legend cell align=right,
legend style={fill=none,draw=none, at={(0.56,1)},
anchor=north,
font=\tiny},
name = ConvD_leg_pos,
font=\scriptsize,
width=0.5\textwidth,
height=0.28\textwidth
]
\addplot [blue, dashed, line width=0.35mm, smooth] table [x index = 0, y index = 1] {GUB_PH_LT_MS_10_4.dat};
\addlegendentry{\text{GUB}, $\rho=10^{4}$}
\addplot [blue, line width=0.35mm, smooth] table [x index = 0, y index = 1] {UB_PH_LT_MS_10_4.dat};
\addlegendentry{\text{UB}, $\rho=10^{4}$}
\end{axis}
\begin{axis}[
hide axis,
legend cell align=right,
legend style={fill=none,draw=none, at={(0.84,1)},
anchor=north,
font=\tiny},
name = ConvD_leg_pos,
font=\scriptsize,
width=0.5\textwidth,
height=0.28\textwidth
]
\addplot [cyan, dashed, line width=0.35mm, smooth] table [x index = 0, y index = 1] {GUB_PH_LT_MS_10_3.dat};
\addlegendentry{\text{GUB}, $\rho=10^{3}$}
\addplot [cyan, line width=0.35mm, smooth] table [x index = 0, y index = 1] {UB_PH_LT_MS_10_3.dat};
\addlegendentry{\text{UB}, $\rho=10^{3}$}
\end{axis}
\begin{axis}[
hide axis,
legend cell align=right,
legend style={fill=none,draw=none, at={(0.79,0.775)},
anchor=north,
font=\tiny},
name = ConvD_leg_pos,
font=\scriptsize,
width=0.5\textwidth,
height=0.28\textwidth
]
\addplot [red, line width=0.1mm, smooth] table [x index = 0, y index = 1] {Direct_opt.dat};
\addlegendentry{Optimal}
\end{axis}
\end{tikzpicture}
}
\caption{Impact of penalty factor $\rho$ on the gap between two bounds on the optimal objective function value. This gap is zero in the optimum.}
\label{fig_gap_iter}
\vspace{-0.5cm}
\end{figure}
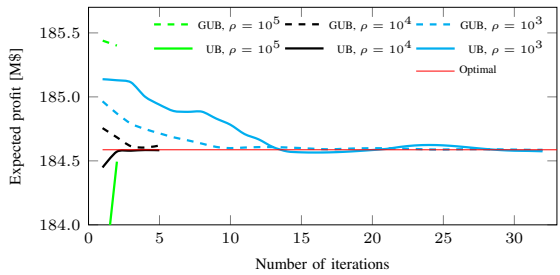

Finally, we show the evolution of both bounds on the optimal objective function value in Fig. \ref{fig_gap_iter}. It shows how the quality of the solution could be traced over iterations depending on the distance between two bounds. At the very first iteration, the nodes of both long- and short-term decisions trees are not tight enough that results in a large gap between the two bounds. This gap reduces over iterations while each scenario-specific investment decision is driven towards the consensus. For small values of penalty factor, both bounds eventually coincide in the global optimum, empirically ensuring the optimality of the solution.

\section{Conclusion} \label{Sec_5}

This paper proposes a suitable consensus-ADMM algorithm to improve the computational tractability of the strategic investment problems in electricity markets. It is based on the relaxation of non-anticipativity conditions of both short- and long-term decision trees of a power producer and their restoration over iterations. Using the proposed algorithm, a decision-maker could include large sets of uncertainties without resorting to restrictive modeling assumptions. \textcolor{blue}{Due to} non-convexity of the original bilevel problem, we introduce a performance guarantee based on the tightness of two bounds on the optimal solution. The algorithm proves to converge to the global optimal solution with a proper tuning of ADMM parameters. Particularly, we show that even with small values of penalty factor, the algorithm results in the optimal solution with the simulation time around 50\% of that provided by the extensive formulation. The algorithm drastically reduces the execution time, \textcolor{blue}{e.g., from 27 minutes to 9 seconds, yielding a near-optimal solution with a relative gap between the two bounds of 0.5\%}.

\section*{Acknowledgment}
The authors would like to thank
A.J. Conejo from the Ohio State University,
T.K. Boomsma from the University of Copenhagen and
T.V. Jensen from the Technical University of Denmark
for thoughtful discussions and constructive criticism of this paper.

\end{document}